# Sparsity oracle inequalities for the Lasso

## Florentina Bunea[*]

*Department of Statistics, Florida State University*
*e-mail:* `bunea@stat.fsu.edu`

## Alexandre Tsybakov

*Laboratoire de Probabilités et Modèles Aléatoires, Université Paris VI*
*e-mail:* `tsybakov@ccr.jussieu.fr`

## Marten Wegkamp[*]

*Department of Statistics, Florida State University, Tallahassee, Florida 32306–4330*
*e-mail:* `wegkamp@stat.fsu.edu`

**Abstract:** This paper studies oracle properties of $\ell_1$-penalized least squares in nonparametric regression setting with random design. We show that the penalized least squares estimator satisfies sparsity oracle inequalities, i.e., bounds in terms of the number of non-zero components of the oracle vector. The results are valid even when the dimension of the model is (much) larger than the sample size and the regression matrix is not positive definite. They can be applied to high-dimensional linear regression, to nonparametric adaptive regression estimation and to the problem of aggregation of arbitrary estimators.

**AMS 2000 subject classifications:** Primary 62G08; secondary 62C20, 62G05, 62G20.

**Keywords and phrases:** sparsity, oracle inequalities, Lasso, penalized least squares, nonparametric regression, dimension reduction, aggregation, mutual coherence, adaptive estimation.

## 1. Introduction

### 1.1. Background

The need for easily implementable methods for regression problems with large number of variables gave rise to an extensive, and growing, literature over the last decade. Penalized least squares with $\ell_1$-type penalties is among the most popular techniques in this area. This method is closely related to restricted least squares minimization, under an $\ell_1$-restriction on the regression coefficients which is called the Lasso method, following [24]. We refer to both methods as Lasso-type methods. Within the linear regression framework these methods became most popular. Let $(Z_1, Y_1), \ldots, (Z_n, Y_n)$ be a sample of independent random

*The research by F. Bunea and M. Wegkamp is supported in part by NSF grant DMS 0706829.

pairs, with $Z_i = (Z_{1i}, \ldots, Z_{Mi})$ and $Y_i = \lambda_1 Z_{1i} + \ldots + \lambda_M Z_{Mi} + W_i$, $\quad i = 1, \cdots, n$, where $W_i$ are independent error terms. Then, for a given $T > 0$, the Lasso estimate of $\lambda \in \mathbb{R}^M$ is

$$\widehat{\lambda}_{\text{lasso}} = \underset{|\lambda|_1 \leq T}{\arg\min} \left\{ \frac{1}{n} \sum_{i=1}^{n} (Y_i - \lambda_1 Z_{1i} - \ldots - \lambda_M Z_{Mi})^2 \right\} \tag{1.1}$$

where $|\lambda|_1 = \sum_{j=1}^{M} |\lambda_j|$. For a given tuning parameter $\gamma > 0$, the penalized estimate of $\lambda \in \mathbb{R}^M$ is

$$\widehat{\lambda}_{\text{pen}} = \underset{\lambda \in \mathbb{R}^M}{\arg\min} \left\{ \frac{1}{n} \sum_{i=1}^{n} (Y_i - \lambda_1 Z_{1i} - \cdots - \lambda_M Z_{Mi})^2 + \gamma |\lambda|_1 \right\}. \tag{1.2}$$

Lasso-type methods can be also applied in the nonparametric regression model $Y = f(X) + W$, where $f$ is the unknown regression function and $W$ is an error term. They can be used to create estimates for $f$ that are linear combinations of basis functions $\phi_1(X), \ldots, \phi_M(X)$ (wavelets, splines, trigonometric polynomials, etc). The vectors of linear coefficients are given by either the $\widehat{\lambda}_{\text{pen}}$ or the $\widehat{\lambda}_{\text{lasso}}$ above, obtained by replacing $Z_{ji}$ by $\phi_j(X_i)$.

In this paper we analyze $\ell_1$-penalized least squares procedures in a more general framework. Let $(X_1, Y_1), \ldots, (X_n, Y_n)$ be a sample of independent random pairs distributed as $(X, Y) \in (\mathcal{X}, \mathbb{R})$, where $\mathcal{X}$ is a Borel subset of $\mathbb{R}^d$; we denote the probability measure of $X$ by $\mu$. Let $f(X) = \mathbb{E}(Y|X)$ be the unknown regression function and $\mathcal{F}_M = \{f_1, \ldots, f_M\}$ be a finite dictionary of real-valued functions $f_j$ that are defined on $\mathcal{X}$. Depending on the statistical targets, the dictionary $\mathcal{F}_M$ can be of different nature. The main examples are:

(I) a collection $\mathcal{F}_M$ of basis functions used to approximate $f$ in the nonparametric regression model as discussed above; these functions need not be orthonormal;

(II) a vector of $M$ one-dimensional random variables $Z = (f_1(X), \ldots, f_M(X))$ as in linear regression;

(III) a collection $\mathcal{F}_M$ of $M$ arbitrary estimators of $f$.

Case (III) corresponds to the aggregation problem: the estimates can arise, for instance, from $M$ different methods; they can also correspond to $M$ different values of the tuning parameter of the same method; or they can be computed on $M$ different data sets generated from the distribution of $(X, Y)$. Without much loss of generality, we treat these estimates $f_j$ as fixed functions; otherwise one can regard our results conditionally on the data set on which they have been obtained.

Within this framework, we use a data dependent $\ell_1$-penalty that differs from the one described in (1.2) in that the tuning parameter $\gamma$ changes with $j$ as in [5, 6]. Formally, for any $\lambda = (\lambda_1, \ldots, \lambda_M) \in \mathbb{R}^M$, define $f_\lambda(x) = \sum_{j=1}^{M} \lambda_j f_j(x)$. Then the penalized least squares estimator of $\lambda$ is

$$\widehat{\lambda} = \underset{\lambda \in \mathbb{R}^M}{\arg\min} \left\{ \frac{1}{n} \sum_{i=1}^{n} \{Y_i - f_\lambda(X_i)\}^2 + \text{pen}(\lambda) \right\}, \tag{1.3}$$

where

$$\text{pen}(\lambda) = 2 \sum_{j=1}^{M} \omega_{n,j} |\lambda_j| \quad \text{with} \quad \omega_{n,j} = r_{n,M} \|f_j\|_n, \tag{1.4}$$

where we write $\|g\|_n^2 = n^{-1} \sum_{i=1}^{n} g^2(X_i)$ for the squared empirical $L_2$ norm of any function $g : \mathcal{X} \to \mathbb{R}$. The corresponding estimate of $f$ is $\widehat{f} = \sum_{j=1}^{M} \widehat{\lambda}_j f_j$. The choice of the tuning sequence $r_{n,M} > 0$ will be discussed in Section 2. Following the terminology used in the machine learning literature (see, e.g., [21]) we call $\widehat{f}$ the aggregate and the optimization procedure $\ell_1$-aggregation.

An attractive feature of $\ell_1$-aggregation is computational feasibility. Because the criterion in (1.3) is convex in $\lambda$, we can use a convex optimization procedure to compute $\widehat{\lambda}$. We refer to [10, 26] for detailed analyzes of these optimization problems and fast algorithms.

Whereas the literature on efficient algorithms is growing very fast, the one on the theoretical aspects of the estimates is still emerging. Most of the existing theoretical results have been derived in the particular cases of either linear or nonparametric regression.

In the linear parametric regression model most results are asymptotic. We refer to [16] for the asymptotic distribution of $\widehat{\lambda}_{\text{pen}}$ in deterministic design regression, when $M$ is fixed and $n \to \infty$. In the same framework, [28, 29] state conditions for subset selection consistency of $\widehat{\lambda}_{\text{pen}}$. For random design Gaussian regression, $M = M(n)$ and possibly larger than $n$, we refer to [20] for consistent variable selection, based on $\widehat{\lambda}_{\text{pen}}$. For similar assumptions on $M$ and $n$, but for random pairs $(Y_i, Z_i)$ that do not necessarily satisfy the linear model assumption, we refer to [12] for the consistency of the risk of $\widehat{\lambda}_{\text{lasso}}$.

The Lasso-type methods have also been extensively used in fixed design non-parametric regression. When the design matrix $\sum_{i=1}^{n} Z_i Z_i'$ is the identity matrix, (1.2) leads to soft thresholding. For soft thresholding in the case of Gaussian errors, the literature dates back to [9]. We refer to [2] for bibliography in the intermediate years and for a discussion of the connections between Lasso-type and thresholding methods, with emphasis on estimation within wavelet bases. For general bases, further results and bibliography we refer to [19]. Under the proper choice of $\gamma$, optimal rates of convergence over Besov spaces, up to logarithmic factors, are obtained. These results apply to the models where the functions $f_j$ are orthonormal with respect to the scalar product induced by the empirical norm. For possible departures from the orthonormality assumption we refer to [5, 6]. These two papers establish finite sample oracle inequalities for the empirical error $\|\widehat{f} - f\|_n^2$ and for the $\ell_1$-loss $|\widehat{\lambda} - \lambda|_1$.

Lasso-type estimators in random design non-parametric regression received very little attention. First results on this subject seem to be [14, 21]. In the aggregation framework described above they established oracle inequalities on the mean risk of $\widehat{f}$, for $\widehat{\lambda}_{\text{lasso}}$ corresponding to $T = 1$ and when $M$ can be larger than $n$. However, this gives an approximation of the oracle risk with the slow rate $\sqrt{(\log M)/n}$, which cannot be improved if $\widehat{\lambda}_{\text{lasso}}$ with fixed $T$ is considered [14, 21]. Oracle inequalities for the empirical error $\|\widehat{f} - f\|_n^2$ and for the $\ell_1$-loss $|\widehat{\lambda} - \lambda|_1$ with faster rates are obtained for $\widehat{\lambda} = \widehat{\lambda}_{\mathrm{pen}}$ in [6] but they are operational only when $M < \sqrt{n}$. The paper [15] studies somewhat different estimators involving the $\ell_1$-norms of the coefficients. For a specific choice of basis functions $f_j$ and with $M < \sqrt{n}$ it proves optimal (up to logarithmic factor) rates of convergence of $\widehat{f}$ on the Besov classes without establishing oracle inequalities. Finally we mention the papers [17, 18, 27] that analyze in the same spirit as we do below the sparsity issue for estimators that differ from $\widehat{\lambda}_{\mathrm{pen}}$ in that the goodness-of-fit term in the minimized criterion cannot be the residual sum of squares.

In the present paper we extend the results of [6] in several ways, in particular, we cover sizes $M$ of the dictionary that can be larger than $n$. To our knowledge, theoretical results for $\widehat{\lambda}_{\mathrm{pen}}$ and the corresponding $\widehat{f}$ when $M$ can be larger than $n$ have not been established for random design in either non-parametric regression or aggregation frameworks. Our considerations are related to a remarkable feature of the $\ell_1$-aggregation: $\widehat{\lambda}_{\mathrm{pen}}$, for an appropriate choice of the tuning sequence $r_{n,M}$, has components exactly equal to zero, thereby realizing subset selection. In contrast, for penalties proportional to $\sum_{j=1}^{M} |\lambda_j|^{\alpha}$, $\alpha > 1$, no estimated coefficients will be set to zero in finite samples; see, e.g. [22] for a discussion. The purpose of this paper is to investigate and quantify when $\ell_1$-aggregation can be used as a dimension reduction technique. We address this by answering the following two questions: "When does $\widehat{\lambda} \in \mathbb{R}^M$, the minimizer of (1.3), behave like an estimate in a dimension that is possibly much lower than $M$?" and "When does the aggregate $\widehat{f}$ behave like a linear approximation of $f$ by a smaller number of functions?" We make these questions precise in the following subsection.

## 1.2. Sparsity and dimension reduction: specific targets

We begin by introducing the following notation. Let

$$M(\lambda) = \sum_{j=1}^{M} I_{\{\lambda_j \neq 0\}} = \mathrm{Card} \ J(\lambda)$$

denote the number of non-zero coordinates of $\lambda$, where $I_{\{\cdot\}}$ denotes the indicator function, and $J(\lambda) = \{j \in \{1, \ldots, M\} : \lambda_j \neq 0\}$. The value $M(\lambda)$ characterizes the *sparsity* of the vector $\lambda$: the smaller $M(\lambda)$, the "sparser" $\lambda$.

To motivate and introduce our notion of sparsity we first consider the simple case of linear regression. The standard assumption used in the literature on linear models is $\mathbb{E}(Y|X) = f(X) = \lambda_0' X$, where $\lambda_0 \in \mathbb{R}^M$ has non-zero coefficients only for $j \in J(\lambda_0)$. Clearly, the $\ell_1$-norm $|\widehat{\lambda}_{\mathrm{OLS}} - \lambda_0|_1$ is of order $M/\sqrt{n}$, in probability, if $\widehat{\lambda}_{\mathrm{OLS}}$ is the ordinary least squares estimator of $\lambda_0$ based on *all* $M$ variables. In contrast, the general results of Theorems 1, 2 and 3 below show that $|\widehat{\lambda} - \lambda_0|_1$ is bounded, up to known constants and logarithms, by $M(\lambda_0)/\sqrt{n}$, for $\widehat{\lambda}$ given by (1.3), if in the penalty term (1.4) we take $r_{n,M} = A\sqrt{(\log M)/n}$. This means that the estimator $\hat{\lambda}$ of the parameter $\lambda_0$ adapts to the sparsity of the problem: its estimation error is smaller when the vector $\lambda_0$ is sparser. In other words, we reduce the effective dimension of the problem from $M$ to $M(\lambda_0)$ without any prior knowledge about the set $J(\lambda_0)$ or the value $M(\lambda_0)$. The improvement is particularly important if $M(\lambda_0) \ll M$ .

Since in general $f$ cannot be represented exactly by a linear combination of the given elements $f_j$ we introduce two ways in which $f$ can be close to such a linear combination. The first one expresses the belief that, for some $\lambda^* \in \mathbb{R}^M$, the squared distance from $f$ to $f_{\lambda^*}$ can be controlled, up to logarithmic factors, by $M(\lambda^*)/n$. We call this "weak sparsity". The second one does not involve $M(\lambda^*)$ and states that, for some $\lambda^* \in \mathbb{R}^M$, the squared distance from $f$ to $f_{\lambda^*}$ can be controlled, up to logarithmic factors, by $n^{-1/2}$. We call this "weak approximation".

We now define weak sparsity. Let $C_f > 0$ be a constant depending only on $f$ and

$$\Lambda = \{\lambda \in \mathbb{R}^M : \|f_\lambda - f\|^2 \le C_f r_{n,M}^2 M(\lambda)\} \tag{1.5}$$

which we refer to as the oracle set $\Lambda$. Here and later we denote by $\|\cdot\|$ the $L_2(\mu)$-norm:

$$\|g\|^2 = \int_{\mathcal{X}} g^2(x)\,\mu(dx)$$

and by $<f, g>$ the corresponding scalar product, for any $f, g \in L_2(\mu)$.

If $\Lambda$ is non-empty, we say that $f$ has the *weak sparsity property relative to the dictionary* $\{f_1, \ldots, f_M\}$. We do not need $\Lambda$ to be a large set: card$(\Lambda) = 1$ would suffice. In fact, under the weak sparsity assumption, our targets are $\lambda^*$ and $f^* = f_{\lambda^*}$, with

$$\lambda^* = \arg\min\left\{\|f_\lambda - f\| : \ \lambda \in \mathbb{R}^M, \ M(\lambda) = k^*\right\}$$

where

$$k^* = \min\{M(\lambda) : \ \lambda \in \Lambda\}$$

is the effective or oracle dimension. All the three quantities, $\lambda^*$, $f^*$ and $k^*$, can be considered as oracles. Weak sparsity can be viewed as a milder version of the *strong sparsity* (or simply *sparsity*) property which commonly means that $f$ admits the exact representation $f = f_{\lambda_0}$ for some $\lambda_0 \in \mathbb{R}^M$, with hopefully small $M(\lambda_0)$.

To illustrate the definition of weak sparsity, we consider the framework (I). Then $\|f_\lambda - f\|$ is the approximation error relative to $f_\lambda$ which can be viewed as a "bias term". For many traditional bases $\{f_j\}$ there exist vectors $\lambda$ with the first $M(\lambda)$ non-zero coefficients and other coefficients zero, such that $\|f_\lambda - f\| \le C(M(\lambda))^{-s}$ for some constant $C > 0$, provided that $f$ is a smooth function with $s$ bounded derivatives. The corresponding variance term is typically of the order $M(\lambda)/n$, so that if $r_{n,M} \sim n^{-1/2}$ the relation $\|f_\lambda - f\|^2 \sim r_{n,M}^2 M(\lambda)$ can be viewed as the bias-variance balance realized for $M(\lambda) \sim n^{\frac{1}{2s+1}}$. We will need to

choose $r_{n,M}$ slightly larger,

$$r_{n,M} \sim \sqrt{\frac{\log M}{n}},$$

but this does not essentially affect the interpretation of $\Lambda$. In this example, the fact that $\Lambda$ is non-void means that there exists $\lambda \in \mathbb{R}^M$ that approximately (up to logarithms) realizes the bias-variance balance or at least undersmoothes $f$ (indeed, we have only an inequality between squared bias and variance in the definition of $\Lambda$). Note that, in general, for instance if $f$ is not smooth, the bias-variance balance can be realized on very bad, even inconsistent, estimators.

We define now another oracle set

$$\Lambda' = \left\{ \lambda \in \mathbb{R}^M : \; \|\mathsf{f}_\lambda - f\|^2 \le C'_f r_{n,M} \right\}.$$

If $\Lambda'$ is non-empty, we say that $f$ has the *weak approximation property relative to the the dictionary* $\{f_1, \ldots, f_M\}$. For instance, in the framework (III) related to aggregation $\Lambda'$ is non-empty if we consider functions $f$ that admit $n^{-1/4}$-consistent estimators in the set of linear combinations $\mathsf{f}_\lambda$, for example, if at least one of the $f_j$'s is $n^{-1/4}$-consistent. This is a modest rate, and such an assumption is quite natural if we work with standard regression estimators $f_j$ and functions $f$ that are not extremely non-smooth.

We will use the notion of weak approximation only in the mutual coherence setting that allows for mild correlation among the $f_j$'s and is considered in Section 2.2 below. Standard assumptions that make our finite sample results work in the asymptotic setting, when $n \to \infty$ and $M \to \infty$, are:

$$r_{n,M} = A\sqrt{\frac{\log M}{n}}$$

for some sufficiently large $A$ and

$$M(\lambda) \le A'\sqrt{\frac{n}{\log M}}$$

for some sufficiently small $A'$, in which case all $\lambda \in \Lambda$ satisfy

$$\|\mathsf{f}_\lambda - f\|^2 \le C'_f r_{n,M}$$

for some constant $C'_f > 0$ depending only on $f$, and weak approximation follows from weak sparsity. However, in general, $r_{n,M}$ and $C_f r_{n,M}^2 M(\lambda)$ are not comparable. So it is not true that weak sparsity implies weak approximation or vice versa. In particular, $C_f r_{n,M}^2 M(\lambda) \le r_{n,M}$, only if $M(\lambda)$ is smaller in order than $\sqrt{n/\log(M)}$, for our choice for $r_{n,M}$.

### *1.3. General assumptions*

We begin by listing and commenting on the assumptions used throughout the paper.

The first assumption refers to the error terms $W_i = Y_i - f(X_i)$. We recall that $f(X) = \mathbb{E}(Y | X)$.

ASSUMPTION (A1). *The random variables $X_1, \ldots, X_n$ are independent, identically distributed random variables with probability measure $\mu$. The random variables $W_i$ are independently distributed with*

$$\mathbb{E}\{W_i \,|\, X_1, \ldots, X_n\} = 0$$

*and*

$$\mathbb{E}\left\{\exp(|W_i|) \,|\, X_1, \ldots, X_n\right\} \leq b \;\; \text{for some finite } b > 0 \text{ and } i = 1, \ldots, n.$$

We also impose mild conditions on $f$ and on the functions $f_j$. Let $\|g\|_\infty = \sup_{x \in \mathcal{X}} |g(x)|$ for any bounded function $g$ on $\mathcal{X}$.

ASSUMPTION (A2). *(a) There exists $0 < L < \infty$ such that $\|f_j\|_\infty \leq L$ for all $1 \leq j \leq M$.*
*(b) There exists $c_0 > 0$ such that $\|f_j\| \geq c_0$ for all $1 \leq j \leq M$.*
*(c) There exists $L_0 < \infty$ such that $\mathbb{E}[f_i^2(X) f_j^2(X)] \leq L_0$ for all $1 \leq i, j \leq M$.*
*(d) There exists $L_* < \infty$ such that $\|f\|_\infty \leq L_* < \infty$.*

*Remark 1.* We note that $(a)$ trivially implies $(c)$. However, as the implied bound may be too large, we opted for stating $(c)$ separately. Note also that $(a)$ and $(d)$ imply the following: for any fixed $\lambda \in \mathbb{R}^M$, there exists a positive constant $L(\lambda)$, depending on $\lambda$, such that $\|f - \mathsf{f}_\lambda\|_\infty = L(\lambda)$.

## 2. Sparsity oracle inequalities

In this section we state our results. They have the form of *sparsity oracle inequalities* that involve the value $M(\lambda)$ in the bounds for the risk of the estimators. All the theorems are valid for arbitrary fixed $n \geq 1, M \geq 2$ and $r_{n,M} > 0$.

### 2.1. Weak sparsity and positive definite inner product matrix

The further analysis of the $\ell_1$-aggregate depends crucially on the behavior of the $M \times M$ matrix $\Psi_M$ given by

$$\Psi_M = \left(\mathbb{E} f_j(X) f_{j'}(X)\right)_{1 \leq j, j' \leq M} = \left(\int f_j(x) f_{j'}(x) \, \mu(dx)\right)_{1 \leq j, j' \leq M}.$$

In this subsection we consider the following assumption

ASSUMPTION (A3). *For any $M \geq 2$ there exist constants $\kappa_M > 0$ such that*

$$\Psi_M - \kappa_M \operatorname{diag}(\Psi_M)$$

*is positive semi-definite.*

Note that $0 < \kappa_M \leq 1$. We will always use Assumption (A3) coupled with (A2). Clearly, Assumption (A3) and part (b) of (A2) imply that the matrix $\Psi_M$ is positive definite, with the minimal eigenvalue $\tau$ bounded from below by $c_0 \kappa_M$. Nevertheless, we prefer to state both assumptions separately, because this allows us to make more transparent the role of the (potentially small) constants $c_0$ and $\kappa_M$ in the bounds, rather than working with $\tau$ which can be as small as their product.

**Theorem 2.1.** *Assume that (A1) – (A3) hold. Then, for all $\lambda \in \Lambda$ we have*

$$\mathbb{P}\left\{\|\widehat{f} - f\|^2 \leq B_1 \kappa_M^{-1} r_{n,M}^2 M(\lambda)\right\} \geq 1 - \pi_{n,M}(\lambda)$$

*and*

$$\mathbb{P}\left\{|\widehat{\lambda} - \lambda|_1 \leq B_2 \kappa_M^{-1} r_{n,M} M(\lambda)\right\} \geq 1 - \pi_{n,M}(\lambda)$$

*where $B_1 > 0$ and $B_2 > 0$ are constants depending on $c_0$ and $C_f$ only and*

$$
\begin{aligned}
\pi_{n,M}(\lambda) &\leq 10M^2 \exp\left(-c_1 n \min\left\{r_{n,M}^2, \frac{r_{n,M}}{L}, \frac{1}{L^2}, \frac{\kappa_M^2}{L_0 M^2(\lambda)}, \frac{\kappa_M}{L^2 M(\lambda)}\right\}\right) \\
&\quad + \exp\left(-c_2 \frac{M(\lambda)}{L^2(\lambda)} n r_{n,M}^2\right),
\end{aligned}
$$

*for some positive constants $c_1, c_2$ depending on $c_0, C_f$ and $b$ only and $L(\lambda) = \|f - \mathsf{f}_\lambda\|_\infty$.*

Since we favored readable results and proofs over optimal constants, not too much attention should be paid to the values of the constants involved. More details about the constants can be found in Section 4.

The most interesting case of Theorem 2.1 corresponds to $\lambda = \lambda^*$ and $M(\lambda) \leq M(\lambda^*) = k^*$. In view of Assumption (A2) we also have a rough bound $L(\lambda^*) \leq L_* + L|\lambda^*|_1$ which can be further improved in several important examples, so that $M(\lambda^*)$ and not $|\lambda^*|_1$ will be involved (cf. Section 3).

### 2.2. Weak sparsity and mutual coherence

The results of the previous subsection hold uniformly over $\lambda \in \Lambda$, when the approximating functions satisfy assumption (A3). We recall that implicit in the definition of $\Lambda$ is the fact that $f$ is well approximated by a smaller number of the given functions $f_j$. Assumption (A3) on the matrix $\Psi_M$ is, however, independent of $f$.

A refinement of our sparsity results can be obtained for $\lambda$ in a set $\Lambda_1$ that combines the requirements for $\Lambda$, while replacing (A3) by a condition on $\Psi_M$ that also depends on $M(\lambda)$. Following the terminology of [8], we consider now matrices $\Psi_M$ with mutual coherence property. We will assume that the correlation

$$\rho_M(i, j) = \frac{<f_i, f_j>}{\|f_i\| \|f_j\|}$$

between elements $i \neq j$ is relatively small, for $i \in J(\lambda)$. Our condition is somewhat weaker than the mutual coherence property defined in [8] where all the correlations for $i \neq j$ are supposed to be small. In our setting the correlations $\rho_M(i,j)$ with $i,j \notin J(\lambda)$ can be arbitrarily close to 1 or to $-1$. Note that such $\rho_M(i,j)$ constitute the overwhelming majority of the elements of the correlation matrix if $J(\lambda)$ is a set of small cardinality: $M(\lambda) \ll M$.

Set

$$\rho(\lambda) = \max_{i \in J(\lambda)} \max_{j \neq i} |\rho_M(i,j)|.$$

With $\Lambda$ given by (1.5) define

$$\Lambda_1 = \{\lambda \in \Lambda : \rho(\lambda)M(\lambda) \leq 1/45\}. \tag{2.1}$$

**Theorem 2.2.** *Assume that (A1) and (A2) hold. Then, for all $\lambda \in \Lambda_1$ we have, with probability at least $1 - \widetilde{\pi}_{n,M}(\lambda)$,*

$$\|\widehat{f} - f\|^2 \leq Cr_{n,M}^2 M(\lambda)$$

*and*

$$|\widehat{\lambda} - \lambda|_1 \leq Cr_{n,M}M(\lambda),$$

*where $C > 0$ is a constant depending only on $c_0$ and $C_f$, and $\widetilde{\pi}_{n,M}(\lambda)$ is defined as $\pi_{n,M}(\lambda)$ in Theorem 2.1 with $\kappa_M = 1$.*

Note that in Theorem 2.2 we do not assume positive definiteness of the matrix $\Psi_M$. However, it is not hard to see that the condition $\rho(\lambda)M(\lambda) \leq 1/45$ implies positive definiteness of the "small" $M(\lambda) \times M(\lambda)$-dimensional submatrix ($< f_i, f_j >)_{i,j \in J(\lambda)}$ of $\Psi_M$.

The numerical constant $1/45$ is not optimal. It can be multiplied at least by a factor close to 4 by taking constant factors close to 1 in the definition of the set $E_2$ in Section 4. The price to pay is a smaller value of constant $c_1$ in the probability $\widetilde{\pi}_{n,M}(\lambda)$.

### 2.3. Weak approximation and mutual coherence

For $\Lambda'$ given in the Introduction, define

$$\Lambda_2 = \{\lambda \in \Lambda' : \rho(\lambda)M(\lambda) \leq 1/45\}. \tag{2.2}$$

**Theorem 2.3.** *Assume that (A1) and (A2) hold. Then, for all $\lambda \in \Lambda_2$, we have*

$$\mathbb{P}\left[\|\widehat{f} - f\|^2 + r_{n,M}|\widehat{\lambda} - \lambda|_1 \leq C'\left\{\|f_\lambda - f\|^2 + r_{n,M}^2 M(\lambda)\right\}\right] \geq 1 - \pi'_{n,M}(\lambda)$$

*where $C' > 0$ is a constant depending only on $c_0$ and $C'_f$, and*

$$\begin{aligned}
\pi'_{n,M}(\lambda) &\leq 14M^2 \exp\left(-c'_1 n \min\left\{\frac{r_{n,M}^2}{L_0}, \frac{r_{n,M}}{L^2}, \frac{1}{L_0 M^2(\lambda)}, \frac{1}{L^2 M(\lambda)}\right\}\right) \\
&\quad + \exp\left(-c'_2 \frac{M(\lambda)}{L^2(\lambda)} n r_{n,M}^2\right)
\end{aligned}$$

*for some constants $c'_1, c'_2$ depending on $c_0, C'_f$ and $b$ only.*

Theorems 2.1 – 2.3 are non-asymptotic results valid for any $r_{n,M} > 0$. If we study asymptotics when $n \to \infty$ or both $n$ and $M$ tend to $\infty$, the optimal choice of $r_{n,M}$ becomes a meaningful question. It is desirable to choose the smallest $r_{n,M}$ such that the probabilities $\pi_{n,M}, \bar{\pi}_{n,M}, \pi'_{n,M}$ tend to 0 (or tend to 0 at a given rate if such a rate is specified in advance). A typical application is in the case where $n \to \infty$, $M = M_n \to \infty$, $\kappa_M$ (when using Theorem 2.1), $L_0, L, L(\lambda^*)$ are independent of $n$ and $M$, and

$$\frac{n}{M^2(\lambda^*) \log M} \to \infty, \quad \text{as } n \to \infty. \tag{2.3}$$

In this case the probabilities $\pi_{n,M}, \bar{\pi}_{n,M}, \pi'_{n,M}$ tend to 0 as $n \to \infty$ if we choose

$$r_{n,M} = A \sqrt{\frac{\log M}{n}}$$

for some sufficiently large $A > 0$. Condition (2.3) is rather mild. It implies, however, that $M$ cannot grow faster than an exponent of $n$ and that $M(\lambda^*) = o(\sqrt{n})$.

## 3. Examples

### 3.1. High-dimensional linear regression

The simplest example of application of our results is in linear parametric regression where the number of covariates $M$ can be much larger than the sample size $n$. In our notation, linear regression corresponds to the case where there exists $\lambda^* \in \mathbb{R}^M$ such that $f = f_{\lambda^*}$. Then the weak sparsity and the weak approximation assumptions hold in an obvious way with $C_f = C'_f = 0$, whereas $L(\lambda^*) = 0$, so that we easily get the following corollary of Theorems 2.1 and 2.2.

**Corollary 1.** *Let $f = f_{\lambda^*}$ for some $\lambda^* \in \mathbb{R}^M$. Assume that (A1) and items (a) – (c) of (A2) hold.*

*(i) If (A3) is satisfied, then*

$$\mathbb{P}\left\{(\widehat{\lambda} - \lambda^*)' \Psi_M (\widehat{\lambda} - \lambda^*) \le B_1 \kappa_M^{-1} r_{n,M}^2 M(\lambda^*)\right\} \ge 1 - \pi^*_{n,M} \tag{3.1}$$

*and*

$$\mathbb{P}\left\{|\widehat{\lambda} - \lambda^*|_1 \le B_2 \kappa_M^{-1} r_{n,M} M(\lambda^*)\right\} \ge 1 - \pi^*_{n,M} \tag{3.2}$$

*where $B_1 > 0$ and $B_2 > 0$ are constants depending on $c_0$ only and*

$$\pi^*_{n,M} \le 10 M^2 \exp\left(-c_1 n \min\left\{r_{n,M}^2, \frac{r_{n,M}}{L}, \frac{1}{L^2}, \frac{\kappa_M^2}{L_0 M^2(\lambda^*)}, \frac{\kappa_M}{L^2 M(\lambda^*)}\right\}\right)$$

*for a positive constant $c_1$ depending on $c_0$ and $b$ only.*

(ii) *If the mutual coherence assumption* $\rho(\lambda^*)M(\lambda^*) \leq 1/45$ *is satisfied, then (3.1) and (3.2) hold with* $\kappa_M = 1$ *and*

$$\pi^*_{n,M} \quad \leq \quad 10M^2 \exp\left(-c_1 n \min\left\{r^2_{n,M}, \frac{r_{n,M}}{L}, \frac{1}{L_0 M^2(\lambda^*)}, \frac{1}{L^2 M(\lambda^*)}\right\}\right)$$

*for a positive constant* $c_1$ *depending on* $c_0$ *and* $b$ *only.*

Result (3.2) can be compared to [7] which gives a control on the $\ell_2$ (not $\ell_1$) deviation between $\widehat{\lambda}$ and $\lambda^*$ in the linear parametric regression setting when $M$ can be larger than $n$, for a different estimator than ours. Our analysis is in several aspects more involved than that in [7] because we treat the regression model with random design and do not assume that the errors $W_i$ are Gaussian. This is reflected in the structure of the probabilities $\pi^*_{n,M}$. For the case of Gaussian errors and fixed design considered in [7], sharper bounds can be obtained (cf. [5]).

### 3.2. Nonparametric regression and orthonormal dictionaries

Assume that the regression function $f$ belongs to a class of functions $\mathcal{F}$ described by some smoothness or other regularity conditions arising in nonparametric estimation. Let $\mathcal{F}_M = \{f_1, \ldots, f_M\}$ be the first $M$ functions of an orthonormal basis $\{f_j\}_{j=1}^{\infty}$. Then $\widehat{f}$ is an estimator of $f$ obtained by an expansion w.r.t. to this basis with data dependent coefficients. Previously known methods of obtaining reasonable estimators of such type for regression with random design mainly have the form of least squares procedures on $\mathcal{F}$ or on a suitable sieve (these methods are not adaptive since $\mathcal{F}$ should be known) or two-stage adaptive procedures where on the first stage least squares estimators are computed on suitable subsets of the dictionary $\mathcal{F}_M$; then, on the second stage, a subset is selected in a data-dependent way, by minimizing a penalized criterion with the penalty proportional to the dimension of the subset. For an overview of these methods in random design regression we refer to [3], to the book [13] and to more recent papers [4, 15] where some other methods are suggested. Note that penalizing by the dimension of the subset as discussed above is not always computationally feasible. In particular, if we need to scan all the subsets of a huge dictionary, or at least all its subsets of large enough size, the computational problem becomes NP-hard. In contrast, the $\ell_1$-penalized procedure that we consider here is computationally feasible. We cover, for example, the case where $\mathcal{F}$'s are the $L_0(\cdot)$ classes (see below). Results of Section 2 imply that an $\ell_1$-penalized procedure is adaptive on the scale of such classes. This can be viewed as an extension to a more realistic random design regression model of Gaussian sequence space results in [1, 11]. However, unlike some results obtained in these papers, we do not establish sharp asymptotics of the risks.

To give precise statements, assume that the distribution $\mu$ of $X$ admits a density w.r.t. the Lebesgue measure which is bounded away from zero by $\mu_{\min} > 0$ and bounded from above by $\mu_{\max} < \infty$. Assume that $\mathcal{F}_M = \{f_1, \ldots, f_M\}$ is an orthonormal system in $L_2(\mathcal{X}, dx)$. Clearly, item (b) of Assumption (A2) holds with $c_0 = \mu_{\min}$, the matrix $\Psi_M$ is positive definite and Assumption (A3) is satisfied with $\kappa_M$ independent of $n$ and $M$. Therefore, we can apply Theorem 2.1. Furthermore, Theorem 2.1 remains valid if we replace there $\|\cdot\|$ by $\|\cdot\|_{\text{Leb}}$ which is the norm in $L_2(\mathcal{X}, dx)$. In this context, it is convenient to redefine the oracle $\lambda^*$ in an equivalent form:

$$\lambda^* = \arg\min \left\{ \|f_\lambda - f\|_{\text{Leb}} : \ \lambda \in \mathbb{R}^M, \ M(\lambda) = k^* \right\} \qquad (3.3)$$

with $k^*$ as before. It is straightforward to see that the oracle (3.3) can be explicitly written as $\lambda^* = (\lambda_1^*, \ldots, \lambda_M^*)$ where $\lambda_j^* = <f_j, f>_{\text{Leb}}$ if $|<f_j, f>_{\text{Leb}}|$ belongs to the set of $k^*$ maximal values among

$$|<f_1, f>_{\text{Leb}}|, \ldots, |<f_M, f>_{\text{Leb}}|$$

and $\lambda_j^* = 0$ otherwise. Here $<\cdot, \cdot>_{\text{Leb}}$ is the scalar product induced by the norm $\|\cdot\|_{\text{Leb}}$. Note also that if $\|f\|_\infty \leq L_*$ we have $L(\lambda^*) = O(M(\lambda^*))$. In fact, $L(\lambda^*) \triangleq \|f - f_{\lambda^*}\| \leq L_* + L|\lambda^*|_1$, whereas

$$
\begin{aligned}
|\lambda^*|_1 &\leq M(\lambda^*) \max_{1 \leq j \leq M} |<f_j, f>_{\text{Leb}}| \leq \frac{M(\lambda^*)}{\mu_{\min}} \max_{1 \leq j \leq M} |<f_j, f>| \\
&\leq \frac{M(\lambda^*) L_* L}{\mu_{\min}}.
\end{aligned}
$$

In the remainder of this section we consider the special case where $\{f_j\}_{j=0}^\infty$ is the Fourier basis in $L_2[0,1]$ defined by $f_1(x) \equiv 1$, $f_{2k}(x) = \sqrt{2}\cos(2\pi k x)$, $f_{2k+1}(x) = \sqrt{2}\sin(2\pi k x)$ for $k = 1, 2, \ldots$, $x \in [0,1]$, and we choose $r_{n,M} = A\sqrt{\frac{\log n}{n}}$. Set for brevity $\theta_j = <f_j, f>_{\text{Leb}}$ and assume that $f$ belongs to the class

$$L_0(k) = \left\{ f : [0,1] \to \mathbb{R} : \ \text{Card} \{ j : \ \theta_j \neq 0 \} \leq k \right\}$$

where $k$ is an unknown integer.

**Corollary 2.** *Let Assumption (A1) and assumptions of this subsection hold. Let $\gamma < 1/2$ be a given number and $M \leq n^s$ for some $s > 0$. Then, for $r_{n,M} = A\sqrt{\frac{\log n}{n}}$ with $A > 0$ large enough, the estimator $\widehat{f}$ satisfies*

$$\sup_{f \in L_0(k)} \mathbb{P} \left\{ \|\widehat{f} - f\|^2 \leq b_1 A^2 \left( \frac{k \log n}{n} \right) \right\} \geq 1 - n^{-b_2}, \quad \forall \ k \leq n^\gamma, \qquad (3.4)$$

*where $b_1 > 0$ is a constant depending on $\mu_{\min}$ and $\mu_{\max}$ only and $b_2 > 0$ is a constant depending also on $A$, $\gamma$ and $s$.*

Proof of this corollary consists in application of Theorem 2.1 with $M(\lambda^*) = k$ and $L(\lambda^*) = 0$ where the oracle $\lambda^*$ is defined in (3.3).

We finally give another corollary of Theorem 2.1 resulting, in particular, in classical nonparametric rates of convergence, up to logarithmic factors. Consider

the class of functions

$$\mathcal{F} = \left\{ f : [0,1] \to \mathbb{R} : \sum_{j=1}^{\infty} |\theta_j| \leq \bar{L} \right\} \tag{3.5}$$

where $\bar{L} > 0$ is a fixed constant. This is a very large class of functions. It contains, for example, all the periodic Hölderian functions on [0,1] and all the Sobolev classes of functions $\mathcal{F}_\beta = \left\{ f : [0,1] \to \mathbb{R} : \sum_{j=1}^{\infty} j^{2\beta} \theta_j^2 \leq Q \right\}$ with smoothness index $\beta > 1/2$ and $Q = Q(\bar{L}) > 0$.

**Corollary 3.** *Let Assumption (A1) and assumptions of this subsection hold. Let $M \leq n^s$ for some $s > 0$. Then, for $r_{n,M} = A\sqrt{\frac{\log n}{n}}$ with $A > 0$ large enough, the estimator $\widehat{f}$ satisfies*

$$\mathbb{P} \left\{ \|\widehat{f} - f\|^2 \leq b_3 \left( \frac{A^2 \log n}{n} \right) M(\lambda^*) \right\} \geq 1 - \pi_n(\lambda^*), \quad \forall \, f \in \mathcal{F}, \tag{3.6}$$

*where $\lambda^*$ is defined in (3.3), $b_3 > 0$ is a constant depending on $\mu_{\min}$ and $\mu_{\max}$ only and*

$$\pi_n(\lambda^*) \leq n^{-b_4} + M^2 \exp(-b_5 n M^{-2}(\lambda^*))$$

*with the constants $b_4 > 0$ and $b_5 > 0$ depending only on $\mu_{\min}$, $\mu_{\max}$, $A$, $\bar{L}$ and $s$.*

This corollary implies, in particular, that the estimator $\widehat{f}$ adapts to unknown smoothness, up to logarithmic factors, simultaneously on the Hölder and Sobolev classes. In fact, it is not hard to see that, for example, when $f \in \mathcal{F}_\beta$ with $\beta > 1/2$ we have $M(\lambda^*) \leq M_n$ where $M_n \sim (n/\log n)^{1/(2\beta+1)}$. Therefore, Corollary 3 implies that $\widehat{f}$ converges to $f$ with rate $(n/\log n)^{-\beta/(2\beta+1)}$, whatever the value $\beta > 1/2$, thus realizing adaptation to the unknown smoothness $\beta$. Similar reasoning works for the Hölder classes.

# 4. Proofs

## *4.1. Proof of Theorem 1*

Throughout this proof $\lambda$ is an arbitrary, fixed element of $\Lambda$ given in (1.5). Recall the notation $f_\lambda = \sum_{j=1}^{M} \lambda_j f_j$. We begin by proving two lemmas. The first one is an elementary consequence of the definition of $\widehat{\lambda}$. Define the random variables

$$V_j = \frac{1}{n} \sum_{i=1}^{n} f_j(X_i) W_i, \quad 1 \leq j \leq M,$$

and the event

$$E_1 = \bigcap_{j=1}^{M} \left\{ 2|V_j| \leq \omega_{n,j} \right\}.$$

**Lemma 1.** *On the event $E_1$, we have for all $n \geq 1$,*

$$\|\widehat{f} - f\|_n^2 + \sum_{j=1}^{M} \omega_{n,j} |\widehat{\lambda}_j - \lambda_j| \quad \leq \quad \|\mathsf{f}_\lambda - f\|_n^2 + 4 \sum_{j \in J(\lambda)} \omega_{n,j} |\widehat{\lambda}_j - \lambda_j|. \, (4.1)$$

*Proof.* We begin as in [19]. By definition, $\widehat{f} = \mathsf{f}_{\widehat{\lambda}}$ satisfies

$$\widehat{S}(\widehat{\lambda}) + \sum_{j=1}^{M} 2\omega_{n,j} |\widehat{\lambda}_j| \quad \leq \quad \widehat{S}(\lambda) + \sum_{j=1}^{M} 2\omega_{n,j} |\lambda_j|$$

for all $\lambda \in \mathbb{R}^M$, which we may rewrite as

$$\|\widehat{f} - f\|_n^2 + \sum_{j=1}^{M} 2\omega_{n,j} |\widehat{\lambda}_j| \leq \|\mathsf{f}_\lambda - f\|_n^2 + \sum_{j=1}^{M} 2\omega_{n,j} |\lambda_j| + \frac{2}{n} \sum_{i=1}^{n} W_i(\widehat{f} - \mathsf{f}_\lambda)(X_i).$$

If $E_1$ holds we have

$$\frac{2}{n} \sum_{i=1}^{n} W_i(\widehat{f} - \mathsf{f}_\lambda)(X_i) \quad = \quad 2 \sum_{j=1}^{M} V_j(\widehat{\lambda}_j - \lambda_j) \leq \sum_{j=1}^{M} \omega_{n,j} |\widehat{\lambda}_j - \lambda_j|$$

and therefore, still on $E_1$,

$$\|\widehat{f} - f\|_n^2 \quad \leq \quad \|\mathsf{f}_\lambda - f\|_n^2 + \sum_{j=1}^{M} \omega_{n,j} |\widehat{\lambda}_j - \lambda_j| + \sum_{j=1}^{M} 2\omega_{n,j} |\lambda_j| - \sum_{j=1}^{M} 2\omega_{n,j} |\widehat{\lambda}_j|.$$

Adding the term $\sum_{j=1}^{M} \omega_{n,j} |\widehat{\lambda}_j - \lambda_j|$ to both sides of this inequality yields further, on $E_1$,

$$\|\widehat{f} - f\|_n^2 + \sum_{j=1}^{M} \omega_{n,j} |\widehat{\lambda}_j - \lambda_j| \leq$$

$$\|\mathsf{f}_\lambda - f\|_n^2 + 2 \sum_{j=1}^{M} \omega_{n,j} |\widehat{\lambda}_j - \lambda_j| + \sum_{j=1}^{M} 2\omega_{n,j} |\lambda_j| - \sum_{j=1}^{M} 2\omega_{n,j} |\widehat{\lambda}_j|.$$

Recall that $J(\lambda)$ denotes the set of indices of the non-zero elements of $\lambda$, and that $M(\lambda) = \mathrm{Card}\, J(\lambda)$. Rewriting the right-hand side of the previous display, we find that, on $E_1$,

$$\|\widehat{f} - f\|_n^2 + \sum_{j=1}^{M} \omega_{n,j} |\widehat{\lambda}_j - \lambda_j|$$

$$\leq \quad \|\mathsf{f}_\lambda - f\|_n^2 + \left( \sum_{j=1}^{M} 2\omega_{n,j} |\widehat{\lambda}_j - \lambda_j| - \sum_{j \notin J(\lambda)} 2\omega_{n,j} |\widehat{\lambda}_j| \right)$$

$$+ \left( - \sum_{j \in J(\lambda)} 2\omega_{n,j} |\widehat{\lambda}_j| + \sum_{j \in J(\lambda)} 2\omega_{n,j} |\lambda_j| \right)$$

$$\leq \quad \|\mathsf{f}_\lambda - f\|_n^2 + 4 \sum_{j \in J(\lambda)} \omega_{n,j} |\widehat{\lambda}_j - \lambda_j|$$

by the triangle inequality and the fact that $\lambda_j = 0$ for $j \notin J(\lambda)$.    □

The following lemma is crucial for the proof of Theorem 1.

**Lemma 2.** *Assume that (A1) – (A3) hold. Define the events*

$$E_2 = \left\{ \frac{1}{2} \|f_j\|^2 \leq \|f_j\|_n^2 \leq 2\|f_j\|^2, \; j = 1, \ldots, M \right\}$$

*and*

$$E_3(\lambda) = \left\{ \|\mathsf{f}_\lambda - f\|_n^2 \leq 2\|\mathsf{f}_\lambda - f\|^2 + r_{n,M}^2 M(\lambda) \right\}.$$

*Then, on the set $E_1 \cap E_2 \cap E_3(\lambda)$, we have*

$$\|\widehat{f} - f\|_n^2 + \frac{c_0 r_{n,M}}{\sqrt{2}} |\widehat{\lambda} - \lambda|_1 \leq \tag{4.2}$$

$$2\|\mathsf{f}_\lambda - f\|^2 + r_{n,M}^2 M(\lambda) + 4r_{n,M} \frac{\sqrt{2M(\lambda)}}{\sqrt{\kappa_M}} \|\widehat{f} - \mathsf{f}_\lambda\|.$$

*Proof.* Observe that assumption (A3) implies that, on the set $E_2$,

$$\sum_{j \in J(\lambda)} \omega_{n,j}^2 |\widehat{\lambda}_j - \lambda_j|^2 \quad \leq \quad \sum_{j=1}^{M} \omega_{n,j}^2 |\widehat{\lambda}_j - \lambda_j|^2$$

$$\leq \quad 2r_{n,M}^2 (\widehat{\lambda} - \lambda)' \text{diag}(\Psi_M)(\widehat{\lambda} - \lambda)$$

$$\leq \quad \frac{2r_{n,M}^2}{\kappa_M} \|\widehat{f} - \mathsf{f}_\lambda\|^2.$$

Applying the Cauchy-Schwarz inequality to the last term on the right hand side of (4.1) and using the inequality above we obtain, on the set $E_1 \cap E_2$,

$$\|\widehat{f} - f\|_n^2 + \sum_{j=1}^{M} \omega_{n,j} |\widehat{\lambda}_j - \lambda_j| \leq \|\mathsf{f}_\lambda - f\|_n^2 + 4r_{n,M} \sqrt{\frac{2M(\lambda)}{\kappa_M}} \|\widehat{f} - \mathsf{f}_\lambda\|.$$

Intersect with $E_3(\lambda)$ and use the fact that $\omega_{n,j} \geq c_0 r_{n,M}/\sqrt{2}$ on $E_2$ to derive the claim.    □

*Proof of Theorem 1.* Recall that $\lambda$ is an arbitrary fixed element of $\Lambda$ given in (1.5). Define the set

$$
\begin{aligned}
U(\lambda) \;=\; & \left\{ \mu \in \mathbb{R}^M : \; \|\mathsf{f}_\mu\| \geq r_{n,M}\sqrt{M(\lambda)} \right\} \cap \\
& \left\{ \mu \in \mathbb{R}^M : \; |\mu|_1 \leq \frac{\sqrt{2}}{c_0}\left( (2C_f+1)r_{n,M}M(\lambda) + 4\sqrt{\frac{2M(\lambda)}{\kappa_M}}\|\mathsf{f}_\mu\| \right) \right\}
\end{aligned}
$$

and the event

$$
E_4(\lambda) = \left\{ \sup_{\mu \in U(\lambda)} \left| \frac{\|\mathsf{f}_\mu\|^2 - \|\mathsf{f}_\mu\|_n^2}{\|\mathsf{f}_\mu\|^2} \right| \leq \frac{1}{2} \right\}.
$$

We prove that the statement of the theorem holds on the event

$$
E(\lambda) := E_1 \cap E_2 \cap E_3(\lambda) \cap E_4(\lambda)
$$

and we bound $\mathbb{P}\left[\{E(\lambda)\}^C\right]$ by $\pi_{n,M}(\lambda)$ in Lemmas 5, 6 and 7 below.

First we observe that, on $E(\lambda) \cap \{\|\widehat{f} - \mathsf{f}_\lambda\| \leq r_{n,M}\sqrt{M(\lambda)}\}$, we immediately obtain, for each $\lambda \in \Lambda$,

$$
\begin{aligned}
\|\widehat{f} - f\| \;\leq\; & \|\mathsf{f}_\lambda - f\| + \|\mathsf{f}_\lambda - \widehat{f}\| \\
\leq\; & \|\mathsf{f}_\lambda - f\| + r_{n,M}\sqrt{M(\lambda)} \\
\leq\; & (1 + C_f^{1/2})r_{n,M}\sqrt{M(\lambda)}
\end{aligned} \tag{4.3}
$$

since $\|\mathsf{f}_\lambda - f\|^2 \leq C_f r_{n,M}^2 M(\lambda)$ for $\lambda \in \Lambda$. Consequently, we find further that, on the same event $E(\lambda) \cap \{\|\widehat{f} - \mathsf{f}_\lambda\| \leq r_{n,M}\sqrt{M(\lambda)}\}$,

$$
\|\widehat{f} - f\|^2 \leq 2(1 + C_f)r_{n,M}^2 M(\lambda) =: C_1 r_{n,M}^2 M(\lambda) \leq C_1 \frac{r_{n,M}^2 M(\lambda)}{\kappa_M},
$$

since $0 < \kappa_M \leq 1$. Also, via (4.2) of Lemma 2 above

$$
|\widehat{\lambda} - \lambda|_1 \leq \frac{1}{c_0}\left\{ 2\sqrt{2}C_f + \sqrt{2} + 8 \right\}\frac{r_{n,M}M(\lambda)}{\sqrt{\kappa_M}} =: C_2 \frac{r_{n,M}M(\lambda)}{\sqrt{\kappa_M}} \leq C_2 \frac{r_{n,M}M(\lambda)}{\kappa_M}.
$$

To finish the proof, we now show that the same conclusions hold on the event $E(\lambda) \cap \left\{ \|\widehat{f} - \mathsf{f}_\lambda\| > r_{n,M}\sqrt{M(\lambda)} \right\}$. Observe that $\widehat{\lambda} - \lambda \in U(\lambda)$ by Lemma 2.

Consequently

$$
\begin{aligned}
\frac{1}{2}\|\widehat{f} - \mathsf{f}_\lambda\|^2 \;\le\; & \|\widehat{f} - \mathsf{f}_\lambda\|_n^2 \\
& \text{(by definition of } E_4(\lambda)) \\
\le\; & 2\|f - \mathsf{f}_\lambda\|_n^2 + 2\|\widehat{f} - f\|_n^2 \\
\le\; & 4\|f - \mathsf{f}_\lambda\|^2 + 2r_{n,M}^2 M(\lambda) + 2\|\widehat{f} - f\|_n^2 \\
& \text{(by definition of } E_3(\lambda)) \\
\le\; & 4\|f - \mathsf{f}_\lambda\|^2 + 2r_{n,M}^2 M(\lambda) + \\
& 2\left\{ 2\|f - \mathsf{f}_\lambda\|^2 + r_{n,M}^2 M(\lambda) + 4r_{n,M}\sqrt{\frac{2M(\lambda)}{\kappa_M}}\|\widehat{f} - \mathsf{f}_\lambda\| \right\} \\
& \text{(by Lemma 2)} \\
\le\; & 8\|f - \mathsf{f}_\lambda\|^2 + 4r_{n,M}^2 M(\lambda) + 4^3 r_{n,M}^2 \frac{2M(\lambda)}{\kappa_M} + \frac{1}{4}\|\widehat{f} - \mathsf{f}_\lambda\|^2
\end{aligned}
\tag{4.4}
$$

using $2xy \le 4x^2 + y^2/4$, with $x = 4r_{n,M}\sqrt{2M(\lambda)/\kappa_M}$ and $y = \|\widehat{f} - \mathsf{f}_\lambda\|$. Hence, on the event $E(\lambda) \cap \{\|\widehat{f} - \mathsf{f}_\lambda\| \ge r_{n,M}M(\lambda)\}$, we have that for each $\lambda \in \Lambda$,

$$
\|\widehat{f} - \mathsf{f}_\lambda\| \le 4\{\sqrt{2C_f} + 6\} r_{n,M}\sqrt{\frac{M(\lambda)}{\kappa_M}}.
\tag{4.5}
$$

This and a reasoning similar to the one used in (4.3) yield

$$
\|\widehat{f} - f\|^2 \le \left\{ (1 + 4\sqrt{2})\sqrt{C_f} + 6 \right\}^2 \kappa_M^{-1} r_{n,M}^2 M(\lambda) =: C_3 \frac{r_{n,M}^2 M(\lambda)}{\kappa_M}.
$$

Also, invoking again Lemma 2 in connection with (4.5) we obtain

$$
|\widehat{\lambda} - \lambda|_1 \le \frac{\sqrt{2}}{c_0} \{2C_f + 1 + 32\sqrt{C_f} + 24\sqrt{2}\} \frac{r_{n,M}M(\lambda)}{\kappa_M} =: C_4 \frac{r_{n,M}M(\lambda)}{\kappa_M}.
$$

Take now $B_1 = C_1 \vee C_3$ and $B_2 = C_2 \vee C_4$ to obtain $\|\widehat{f} - f\|^2 \le B_1 \kappa_M^{-1} r_{n,M}^2 M(\lambda)$ and $|\widehat{\lambda} - \lambda|_1 \le B_2 \kappa_M^{-1} r_{n,M} M(\lambda)$. The conclusion of the theorem follows from the bounds on the probabilities of the complements of the events $E_1, E_2, E_3(\lambda)$ and $E_4(\lambda)$ as proved in Lemmas 4, 5, 6 and 7 below. □

The following results will make repeated use of a version of Bernstein's inequality which we state here for ease of reference.

**Lemma 3** (Bernstein's inequality). *Let $\zeta_1, \ldots, \zeta_n$ be independent random variables such that*

$$
\frac{1}{n}\sum_{i=1}^n \mathbb{E}|\zeta_i|^m \le \frac{m!}{2} w^2 d^{m-2}
$$

*for some positive constants $w$ and $d$ and for all integers $m \geq 2$. Then, for any $\varepsilon > 0$ we have*

$$\mathbb{P} \left\{ \sum_{i=1}^{n} (\zeta_i - \mathbb{E}\zeta_i) \geq n\varepsilon \right\} \leq \exp \left( -\frac{n\varepsilon^2}{2(w^2 + d\varepsilon)} \right). \qquad (4.6)$$

**Lemma 4.** *Assume that (A1) and (A2) hold Then, for all $n \geq 1$, $M \geq 2$,*

$$\mathbb{P} \left( E_2^C \right) \leq 2M \exp \left( -\frac{nc_0^2}{12L^2} \right). \qquad (4.7)$$

*Proof.* The proof follows from a simple application of the union bound and Bernstein's inequality:

$$
\begin{aligned}
\mathbb{P} \left( E_2^C \right) &\leq M \max_{1 \leq j \leq M} \left( \mathbb{P} \left\{ \frac{1}{2} \|f_j\|^2 > \|f_j\|_n^2 \right\} + \mathbb{P} \left\{ \|f_j\|_n^2 > 2\|f_j\|^2 \right\} \right) \\
&\leq M \exp \left( -\frac{nc_0^2}{12L^2} \right) + M \exp \left( -\frac{nc_0^2}{4L^2} \right),
\end{aligned}
$$

where we applied Bernstein's inequality with $w^2 = \|f_j\|^2 L^2$ and $d = L^2$ and with $\varepsilon = \frac{1}{2}\|f_j\|^2$ for the first probability and with $\varepsilon = \|f_j\|^2$ for the second one. $\qquad \square$

**Lemma 5.** *Let Assumptions (A1) and (A2) hold. Then*

$$
\begin{aligned}
\mathbb{P} \left( \{E_1 \cap E_2\}^C \right) &\leq 2M \exp \left( -\frac{nr_{n,M}^2}{16b} \right) + 2M \exp \left( -\frac{nr_{n,M}c_0}{8\sqrt{2}L} \right) \\
&\quad + 2M \exp \left( -\frac{nc_0^2}{12L^2} \right).
\end{aligned}
$$

*Proof.* We apply Bernstein's inequality with the variables $\zeta_i = \zeta_{i,j} = f_j(X_i)W_i$, for each fixed $j \in \{1, \dots, M\}$ and fixed $X_1, \dots, X_n$. By assumptions (A1) and (A2), we find that, for $m \geq 2$,

$$
\begin{aligned}
\frac{1}{n} \sum_{i=1}^{n} \mathbb{E} \left\{ |\zeta_{i,j}|^m \mid X_1, \dots, X_n \right\} &\leq L^{m-2} \frac{1}{n} \sum_{i=1}^{n} f_j^2(X_i) \mathbb{E} \left\{ |W_i|^m \mid X_1, \dots, X_n \right\} \\
&\leq \frac{m!}{2} L^{m-2} \left( b \|f_j\|_n^2 \right).
\end{aligned}
$$

Using (4.6), with $\varepsilon = \omega_{n,j}/2$, $w = \sqrt{b}\|f_j\|_n$, $d = L$, the union bound and the fact that

$$\exp\{-x/(\alpha + \beta)\} \leq \exp\{-x/(2\alpha)\} + \exp\{-x/(2\beta)\}, \quad \forall x, \alpha, \beta > 0, \qquad (4.8)$$

we obtain

$$\mathbb{P}\left(E_1^C \mid X_1, \ldots, X_n\right) \leq 2 \sum_{j=1}^{M} \exp\left(-\frac{nr_{n,M}^2 \|f_j\|_n^2/4}{2\left(b\|f_j\|_n^2 + Lr_{n,M}\|f_j\|_n/2\right)}\right)$$

$$\leq 2M \exp\left(-\frac{nr_{n,M}^2}{16b}\right) + 2 \sum_{j=1}^{M} \exp\left(-\frac{nr_{n,M}\|f_j\|_n}{8L}\right).$$

This inequality, together with the fact that on $E_2$ we have $\|f_j\|_n \geq \|f_j\|/\sqrt{2} \geq c_0/\sqrt{2}$, implies

$$\mathbb{P}\left(E_1^C \cap E_2\right) \leq 2M \exp\left(-\frac{nr_{n,M}^2}{16b}\right) + 2M \exp\left(-\frac{nr_{n,M}c_0}{8\sqrt{2}L}\right).$$

Combining this with Lemma 4 we get the result. □

**Lemma 6.** *Assume that (A1) and (A2) hold. Then, for all $n \geq 1$, $M \geq 2$,*

$$\mathbb{P}\left[\{E_3(\lambda)\}^C\right] \leq \exp\left(-\frac{M(\lambda)nr_{n,M}^2}{4L^2(\lambda)}\right).$$

*Proof.* Recall that $\|\mathsf{f}_\lambda - f\|_\infty = L(\lambda)$. The claim follows from Bernstein's inequality applied with $\varepsilon = \|\mathsf{f}_\lambda - f\|^2 + r_{n,M}^2 M(\lambda)$, $d = L^2(\lambda)$ and $w^2 = \|\mathsf{f}_\lambda - f\|^2 L^2(\lambda)$. □

**Lemma 7.** *Assume (A1) – (A3). Then*

$$\mathbb{P}\left[\{E_4(\lambda)\}^C\right] \leq 2M^2 \exp\left(-\frac{n}{16L_0C^2M^2(\lambda)}\right) + 2M^2 \exp\left(-\frac{n}{8L^2CM(\lambda)}\right),$$

*where $C = 2c_0^{-2}\left(2C_f + 1 + 4\sqrt{2/\kappa_M}\right)^2$.*

*Proof.* Let

$$\psi_M(i,j) = \mathbb{E}[f_i(X)f_j(X)] \quad \text{and} \quad \psi_{n,M}(i,j) = \frac{1}{n}\sum_{k=1}^{n} f_i(X_k)f_j(X_k)$$

denote the $(i,j)$th entries of matrices $\Psi_M$ and $\Psi_{n,M}$, respectively. Define

$$\eta_{n,M} = \max_{1 \leq i,j, \leq M} |\psi_M(i,j) - \psi_{n,M}(i,j)|.$$

Then, for every $\mu \in U(\lambda)$ we have

$$
\begin{aligned}
\frac{\left|\|f_\mu\|^2 - \|f_\mu\|_n^2\right|}{\|f_\mu\|^2} &= \frac{|\mu'(\Psi_M - \Psi_{n,M})\mu|}{\|f_\mu\|^2} \\
&\leq \frac{|\mu|_1^2}{\|f_\mu\|^2} \max_{1 \leq i,j, \leq M} |\psi_M(i,j) - \psi_{n,M}(i,j)| \\
&\leq \frac{2}{\|f_\mu\|^2 c_0^2} \left\{ (2C_f + 1)r_{n,M}M(\lambda) + 4\sqrt{\frac{2M(\lambda)}{\kappa_M}}\|f_\mu\| \right\}^2 \eta_{n,M} \\
&\leq \frac{2}{c_0^2} \left\{ (2C_f + 1)\sqrt{M(\lambda)} + 4\sqrt{\frac{2M(\lambda)}{\kappa_M}} \right\}^2 \eta_{n,M} \\
&= \frac{2}{c_0^2} \left\{ (2C_f + 1) + 4\sqrt{\frac{2}{\kappa_M}} \right\}^2 M(\lambda)\eta_{n,M} \\
&= CM(\lambda)\eta_{n,M}
\end{aligned}
$$

Using the the last display and the union bound, we find for each $\lambda \in \Lambda$ that

$$
\begin{aligned}
\mathbb{P}\left[\{E_4(\lambda)\}^C\right] &\leq \mathbb{P}\left[\eta_{n,M} \geq 1/\{2CM(\lambda)\}\right] \\
&\leq 2M^2 \max_{1 \leq i,j \leq M} \mathbb{P}\left[|\psi_M(i,j) - \psi_{n,M}(i,j)| \geq 1/\{2CM(\lambda)\}\right].
\end{aligned}
$$

Now for each $(i,j)$, the value $\psi_M(i,j) - \psi_{n,M}(i,j)$ is a sum of $n$ i.i.d. zero mean random variables. We can therefore apply Bernstein's inequality with $\zeta_k = f_i(X_k)f_j(X_k)$, $\varepsilon = 1/\{2CM(\lambda)\}$, $w^2 = L_0$, $d = L^2$ and inequality (4.8) to obtain the result. $\qquad\square$

### *4.2. Proof of Theorem 2.2*

Let $\lambda$ be an arbitrary fixed element of $\Lambda_1$ given in (2.1). The proof of this theorem is similar to that of Theorem 1. The only difference is that we now show that the result holds on the event

$$
\tilde{E}(\lambda) := E_1 \cap E_2 \cap E_3(\lambda) \cap \tilde{E}_4(\lambda).
$$

Here the set $\tilde{E}_4(\lambda)$ is given by

$$
\widetilde{E}_4(\lambda) = \left\{ \sup_{\mu \in \widetilde{U}(\lambda)} \left| \frac{\|f_\mu\|^2 - \|f_\mu\|_n^2}{\|f_\mu\|^2} \right| \leq \frac{1}{2} \right\},
$$

where

$$
\begin{aligned}
\widetilde{U}(\lambda) &= \left\{ \mu \in \mathbb{R}^M : \|f_\mu\| \geq r_{n,M}\sqrt{M(\lambda)} \right\} \cap \\
&\quad \left\{ \mu \in \mathbb{R}^M : |\mu|_1 \leq \frac{2\sqrt{2}}{c_0} \left( (2C_f + 1)r_{n,M}M(\lambda) + 8\sqrt{M(\lambda)}\|f_\mu\| \right) \right\}.
\end{aligned}
$$

We bounded $\mathbb{P}\left(\{E_1 \cap E_2\}^C\right)$ and $\mathbb{P}\left[\{E_3(\lambda)\}^C\right]$ in Lemmas 5 and 6 above. The bound for $P\left[\{\widetilde{E}_4(\lambda)\}^C\right]$ is obtained exactly as in Lemma 7 but now with $C_1 = 8c_0^{-2}\left(2C_f + 9\right)^2$, so that we have

$$\mathbb{P}\left[\{\widetilde{E}_4(\lambda)\}^C\right] \leq 2M^2 \exp\left(-\frac{n}{16L_0 C_1^2 M^2(\lambda)}\right) + 2M^2 \exp\left(-\frac{n}{8L^2 C_1 M(\lambda)}\right).$$

The proof of Theorem 2.2 on the set $\widetilde{E}(\lambda) \cap \left\{\|\widehat{f} - \mathsf{f}_\lambda\| \leq r_{n,M}\sqrt{M(\lambda)}\right\}$ is identical to that of Theorem 2.1 on the set $E(\lambda) \cap \left\{\|\widehat{f} - \mathsf{f}_\lambda\| \leq r_{n,M}\sqrt{M(\lambda)}\right\}$. Next, on the set $\widetilde{E}(\lambda) \cap \left\{\|\widehat{f} - \mathsf{f}_\lambda\| > r_{n,M}\sqrt{M(\lambda)}\right\}$, we follow again the argument of Theorem 2.1 first invoking Lemma 8 given below to argue that $\widehat{\lambda} - \lambda \in \widetilde{U}(\lambda)$ (this lemma plays the same role as Lemma 2 in the proof of Theorem 2.1) and then reasoning exactly as in (4.4). $\qquad\square$

**Lemma 8.** *Assume that (A1) and (A2) hold and that $\lambda$ is an arbitrary fixed element of the set $\{\lambda \in \mathbb{R}^M : \rho(\lambda)M(\lambda) \leq 1/45\}$. Then, on the set $E_1 \cap E_2 \cap E_3(\lambda)$, we have*

$$\|\widehat{f} - f\|_n^2 + \frac{c_0 r_{n,M}}{2\sqrt{2}}|\widehat{\lambda} - \lambda|_1 \leq \tag{4.9}$$
$$2\|\mathsf{f}_\lambda - f\|^2 + r_{n,M}^2 M(\lambda) + 8r_{n,M}\sqrt{M(\lambda)}\|\widehat{f} - \mathsf{f}_\lambda\|.$$

*Proof.* Set for brevity

$$\rho = \rho(\lambda), \quad u_j = \widehat{\lambda}_j - \lambda_j, \quad a = \sum_{j=1}^{M}\|f_j\|\,|u_j|, \quad a(\lambda) = \sum_{j \in J(\lambda)}\|f_j\|\,|u_j|.$$

By Lemma 1, on $E_1$ we have

$$\|\widehat{f} - f\|_n^2 + \sum_{j=1}^{M}\omega_{n,j}|u_j| \leq \|\mathsf{f}_\lambda - f\|_n^2 + 4\sum_{j \in J(\lambda)}\omega_{n,j}|u_j|. \tag{4.10}$$

Now, on the set $E_2$,

$$4\sum_{j \in J(\lambda)}\omega_{n,j}|u_j| \leq 8r_{n,M}a(\lambda) \leq 8r_{n,M}\sqrt{M(\lambda)}\sqrt{\sum_{j \in J(\lambda)}\|f_j\|^2 u_j^2}. \tag{4.11}$$

Here

$$\begin{aligned}
\sum_{j \in J(\lambda)}\|f_j\|^2 u_j^2 &= \|\widehat{f} - \mathsf{f}_\lambda\|^2 - \sum_{i,j \notin J(\lambda)}<f_i, f_j> u_i u_j \\
&\quad -2\sum_{i \notin J(\lambda)}\sum_{j \in J(\lambda)}<f_i, f_j> u_i u_j - \sum_{i,j \in J(\lambda), i \neq j}<f_i, f_j> u_i u_j \\
&\leq \|\widehat{f} - \mathsf{f}_\lambda\|^2 + 2\rho\sum_{i \notin J(\lambda)}\|f_i\|\,|u_i|\sum_{j \in J(\lambda)}\|f_j\|\,|u_j| + \rho a^2(\lambda) \\
&= \|\widehat{f} - \mathsf{f}_\lambda\|^2 + 2\rho a(\lambda)a - \rho a^2(\lambda)
\end{aligned}$$

where we used the fact that $\sum \sum_{i,j \notin J(\lambda)} < f_i, f_j > u_i u_j \geq 0$. Combining this with the second inequality in (4.11) yields

$$a^2(\lambda) \leq M(\lambda) \left\{ \|\widehat{f} - \mathsf{f}_\lambda\|^2 + 2\rho a(\lambda)a - \rho a^2(\lambda) \right\}$$

which implies

$$a(\lambda) \leq \frac{2\rho M(\lambda)a}{1 + \rho M(\lambda)} + \frac{\sqrt{M(\lambda)}\|\widehat{f} - \mathsf{f}_\lambda\|}{1 + \rho M(\lambda)}. \tag{4.12}$$

¿From (4.10), (4.12) and the first inequality in (4.11) we get

$$\begin{aligned}
\|\widehat{f} - f\|_n^2 + \sum_{j=1}^{M} \omega_{n,j}|u_j| &\leq \|\mathsf{f}_\lambda - f\|_n^2 + \frac{16\rho M(\lambda)r_{n,M}a}{1 + \rho M(\lambda)} \\
&\quad + \frac{8r_{n,M}\sqrt{M(\lambda)}\|\widehat{f} - \mathsf{f}_\lambda\|}{1 + \rho M(\lambda)}.
\end{aligned}$$

Combining this with the fact that $r_{n,M}\|f_j\| \leq \sqrt{2}\omega_{n,j}$ on $E_2$ and $\rho M(\lambda) \leq 1/45$ we find

$$\|\widehat{f} - f\|_n^2 + \frac{1}{2}\sum_{j=1}^{M} \omega_{n,j}|u_j| \quad \leq \quad \|\mathsf{f}_\lambda - f\|_n^2 + 8r_{n,M}\sqrt{M(\lambda)}\|\widehat{f} - \mathsf{f}_\lambda\|.$$

Intersect with $E_3(\lambda)$ and use the fact that $\omega_{n,j} \geq c_0 r_{n,M}/\sqrt{2}$ on $E_2$ to derive the claim. $\qquad\square$

### *4.3. Proof of Theorem 2.3.*

Let $\lambda \in \Lambda_2$ be arbitrary, fixed and we set for brevity $C'_f = 1$. We consider separately the cases (a) $\|\mathsf{f}_\lambda - f\|^2 \leq r_{n,M}^2 M(\lambda)$ and (b) $\|\mathsf{f}_\lambda - f\|^2 > r_{n,M}^2 M(\lambda)$.

*Case (a).* It follows from Theorem 2.2 that

$$\|\widehat{f} - f\|^2 + r_{n,M}|\widehat{\lambda} - \lambda|_1 \leq Cr_{n,M}^2 M(\lambda) < C \left\{ r_{n,M}^2 M(\lambda) + \|\mathsf{f}_\lambda - f\|^2 \right\}$$

with probability greater than $1 - \widetilde{\pi}_{n,M}(\lambda)$.

*Case (b).* In this case it is sufficient to show that

$$\|\widehat{f} - f\|^2 + r_{n,M}|\widehat{\lambda} - \lambda|_1 \leq C'\|\mathsf{f}_\lambda - f\|^2, \tag{4.13}$$

for a constant $C' > 0$, on some event $E'(\lambda)$ with $\mathbb{P}\{E'(\lambda)\} \geq 1 - \pi'_{n,M}(\lambda)$. We proceed as follows. Define the set

$$\begin{aligned}
U'(\lambda) \quad = \quad &\Big\{ \mu \in \mathbb{R}^M : \|\mathsf{f}_\mu\| > \|\mathsf{f}_\lambda - f\|, \\
&|\mu|_1 \leq \frac{2\sqrt{2}}{c_0 r_{n,M}} \left( 3\|\mathsf{f}_\lambda - f\|^2 + 8\|\mathsf{f}_\lambda - f\| \cdot \|\mathsf{f}_\mu\| \right) \Big\}
\end{aligned}$$

and the event

$$E'(\lambda) := E_1 \cap E_2 \cap E_3(\lambda) \cap E_5(\lambda),$$

where

$$E_5(\lambda) = \left\{ \sup_{\mu \in U'(\lambda)} \left| \frac{\|\mathsf{f}_\mu\|^2 - \|\mathsf{f}_\mu\|_n^2}{\|\mathsf{f}_\mu\|^2} \right| \le \frac{1}{2} \right\}.$$

We prove the result by considering two cases separately: $\|\widehat{f} - \mathsf{f}_\lambda\| \le \|\mathsf{f}_\lambda - f\|$ and $\|\widehat{f} - \mathsf{f}_\lambda\| > \|\mathsf{f}_\lambda - f\|$.

On the event $\{\|\widehat{f} - \mathsf{f}_\lambda\| \le \|\mathsf{f}_\lambda - f\|\}$ we have immediately

$$\|\widehat{f} - f\|^2 \le 2\|\widehat{f} - \mathsf{f}_\lambda\|^2 + 2\|\mathsf{f}_\lambda - f\|^2 \le 4\|\mathsf{f}_\lambda - f\|^2. \tag{4.14}$$

Recall that being in *Case (b)* means that $\|\mathsf{f}_\lambda - f\|^2 > r_{n,M}^2 M(\lambda)$. This coupled with (4.14) and with the inequality $\|\widehat{f} - \mathsf{f}_\lambda\| \le \|\mathsf{f}_\lambda - f\|$ shows that the right hand side of (4.9) in Lemma 8 can be bounded, up to multiplicative constants, by $\|\mathsf{f}_\lambda - f\|^2$. Thus, on the event $E'(\lambda) \cap \left\{ \|\widehat{f} - \mathsf{f}_\lambda\| \le \|\mathsf{f}_\lambda - f\| \right\}$ we have

$$r_{n,M}|\widehat{\lambda} - \lambda|_1 \le C\|\mathsf{f}_\lambda - f\|^2,$$

for some constant $C > 0$. Combining this with (4.14) we get (4.13), as desired.

Let now $\|\widehat{f} - \mathsf{f}_\lambda\| > \|\mathsf{f}_\lambda - f\|$. Then, by Lemma 8, we get that $\widehat{\lambda} - \lambda \in U'(\lambda)$, on $E_1 \cap E_2 \cap E_3(\lambda)$. Using this fact and the definition of $E_5(\lambda)$, we find that on $E'(\lambda) \cap \left\{ \|\widehat{f} - \mathsf{f}_\lambda\| > \|\mathsf{f}_\lambda - f\| \right\}$ we have

$$\frac{1}{2}\|\widehat{f} - \mathsf{f}_\lambda\|^2 \le \|\widehat{f} - \mathsf{f}_\lambda\|_n^2.$$

Repeating the argument in (4.4) with the only difference that we use now Lemma 8 instead of Lemma 2 and recalling that $\|\mathsf{f}_\lambda - f\|^2 > r_{n,M}^2 M(\lambda)$ since we are in *Case (b)*, we get

$$\|\widehat{f} - \mathsf{f}_\lambda\|^2 \le C(r_{n,M}^2 M(\lambda) + \|\mathsf{f}_\lambda - f\|^2) \le C''\|\mathsf{f}_\lambda - f\|^2 \tag{4.15}$$

for some constants $C > 0, C'' > 0$. Therefore,

$$\|\widehat{f} - f\|^2 \le 2\|\widehat{f} - \mathsf{f}_\lambda\|^2 + 2\|\mathsf{f}_\lambda - f\|^2 \le (2C'' + 1)\|\mathsf{f}_\lambda - f\|^2. \tag{4.16}$$

Note that (4.15) and (4.16) have the same form (up to multiplicative constants) as the condition $\|\widehat{f} - \mathsf{f}_\lambda\| \le \|\mathsf{f}_\lambda - f\|$ and the inequality (4.14) respectively. Hence, we can use the reasoning following (4.14) to conclude that on $E'(\lambda) \cap \left\{ \|\widehat{f} - \mathsf{f}_\lambda\| > \|\mathsf{f}_\lambda - f\| \right\}$ inequality (4.13) holds true.

The result of the theorem follows now from the bound $\mathbb{P}\left[\{E'(\lambda)\}^C\right] \le \pi'_{n,M}(\lambda)$ which is a consequence of Lemmas 5, 6 and of the next Lemma 9. □

**Lemma 9.** *Assume (A1) and (A2). Then, for all $n \geq 1$, $M \geq 2$,*

$$\mathbb{P}\left[\{E_5(\lambda)\}^C\right] \leq 2M^2 \exp\left(-\frac{nr_{n,M}^2}{16C^2L_0}\right) + 2M^2 \exp\left(-\frac{nr_{n,M}}{8L^2C}\right)$$

*where $C = 8 \cdot 11^2 c_0^{-2}$.*

*Proof.* The proof closely follows that of Lemma 7. Using the inequality $\|\mathsf{f}_\lambda - f\|^2 \leq r_{n,M}$, we deduce that

$$
\begin{aligned}
\mathbb{P}\left[\{E_5(\lambda)\}^C\right] &\leq \mathbb{P}\left\{\eta_{n,M}\frac{8 \cdot 11^2}{c_0^2 r_{n,M}^2}\|\mathsf{f}_\lambda - f\|^2 \geq \frac{1}{2}\right\} \\
&\leq \mathbb{P}\left\{\eta_{n,M}\frac{8 \cdot 11^2}{c_0^2 r_{n,M}} \geq \frac{1}{2}\right\}.
\end{aligned}
$$

An application of Bernstein's inequality with $\zeta_k = f_i(X_k)f_j(X_k)$, $\varepsilon = r_{n,M}/(2C)$, $w^2 = L_0$ and $d = L^2$ completes the proof of the lemma. $\square$